\documentclass[11pt, oneside, reqno]{article}

\usepackage{amssymb}
\usepackage{amsmath}
\usepackage{amsthm}
\usepackage{amsfonts}
\usepackage{enumerate}
\usepackage{esint,subfig}
\usepackage{latexsym,amssymb}
\usepackage{color}
\usepackage[hidelinks]{hyperref}
\usepackage{bbm}
\usepackage[latin1]{inputenc}
\usepackage{multicol}
\usepackage[pdftex]{graphicx}
\usepackage[nottoc,notlot,notlof]{tocbibind}

\usepackage{geometry,color,mathtools,fixltx2e,microtype,lastpage,mathrsfs}

\usepackage{xr}

\allowdisplaybreaks

\numberwithin{equation}{section}

\theoremstyle{plain}
\newtheorem{thm}{Theorem}[section]
\newtheorem{cor}[thm]{Corollary}

\theoremstyle{definition}

\theoremstyle{remark}

\def\R{\mathbb{R}}
\def\P{\mathbb{P}}

\def\S{\mathbb{S}}

\def\EE{\mathcal{E} }
\def\M{\mathscr{M} }
\def\MM{\mathcal{M} }

\def\L{\mathscr{L} }
\def\LL{\mathcal{L} }

\def\E{\mathbb{E}}

\def\to{\longrightarrow}
\def\sto{\rightarrow}

\def\1{\mathbbm{1}}

\DeclareMathOperator{\Cov}{Cov}
\DeclareMathOperator{\Corr}{Corr}
\DeclareMathOperator{\Var}{Var}

\def\area{\operatorname{area}}
\def\diam{\operatorname{diam}}

\DeclareMathOperator{\cum}{cum}

\def\r{\right\rangle}
\def\l{\left\langle}

\makeatletter
\newcommand{\vast}{\bBigg@{3.5}}
\newcommand{\Vast}{\bBigg@{5}}
\makeatother

\renewcommand{\(}{\left(}
\renewcommand{\)}{\right)}

\renewcommand{\]}{\right]}
\renewcommand{\{}{\left\lbrace}
\renewcommand{\}}{\right\rbrace}

\newcommand{\norm}[1]{\left\lVert#1\right\rVert}
\newcommand{\abs}[1]{\left\lvert#1\right\rvert}

\def\endmproof{\hfill \mbox{\raggedright \rule{0.1in}{0.1in}}}

\begin{document}

\title{A Note on the Reduction Principle \\for the Nodal Length of Planar Random Waves}

\date{\footnotesize \today}

\author{{Anna Vidotto}
\smallskip
\\
\smallskip
{\footnotesize \it Universit\`a degli Studi ``G. D'Annunzio'' Chieti-Pescara} \\
{\footnotesize  email: anna.vidotto@unich.it}}

\maketitle

\begin{abstract}
Inspired by the recent work \cite{MRW:19}, we prove that the nodal length of a planar random wave $B_{E}$, i.e. the length of its zero set $B_{E}^{-1}(0)$, is asymptotically equivalent, in the $L^{2}$-sense and in the high-frequency limit $E\sto \infty$, to the integral of $H_{4}(B_{E}(x))$, $H_4$ being the fourth Hermite polynomial. As a straightforward consequence, we obtain a central limit theorem in Wasserstein distance. This complements recent findings in \cite{NPR:19} and \cite{PV:20}.
\end{abstract}

\begin{itemize}
\item \textbf{AMS Classification}: 60G60, 60F05, 34L20, 33C10.

\item \textbf{Keywords and Phrases}: Nodal Length, Random Plane Waves,
Sample Trispectrum, Berry's Cancellation, Quantitative Central Limit Theorem.
\end{itemize}

\section{Introduction and Main Results}

\subsection{Motivation}

Let $(\MM,g)$ be a smooth Riemannian manifold and let $f_k:\MM\to \R$ be a random function which almost surely solves the Helmholtz equation, that is 
$$
\Delta_g f_k+\lambda_kf_k=0 \qquad \text{a.s.}\,,
$$
where $\Delta_g$ is the Laplacian defined with respect to the Riemannian metric $g$ and $-\lambda_k$ its eigenvalue.
The study of the geometric properties of the excursion sets of $f_k$ at a fixed level $u\in \R$, i.e.
$$
\EE_u(f_k,\MM):=\{x\in\MM: f_k(x)\ge u\}\,,
$$
in the high-energy limit $k\sto\infty$, has recently attracted great interest, starting from the seminal work \cite{Be:77}, in which Berry conjectured that, as $k\sto\infty$, local geometric functionals of a planar random eigenfunction $f_k$ reproduce the behavior of a \emph{typical} deterministic Laplace eigenfunction on any \emph{generic} manifold. 
In two dimensions, three important geometric quantities that characterize local geometric functionals associated with a random field are the Euler-Poincar\'e characteristic $\LL_0^{f_k}(\EE_u(f_k,\MM))$, the boundary length $\LL_1^{f_k}(\EE_u(f_k,\MM))$ and the area $\LL_2^{f_k}(\EE_u(f_k,\MM))$, namely, the so-called Lipschitz-Killing curvatures (see \cite{AT:07}). 

Among these geometric functionals, particular attention was drawn by the behavior of the nodal length (the boundary length at $u=0$), starting from the celebrated Yau's conjecture on its value for deterministic eigenfunctions on general manifolds, see \cite{Ya:82}. With a physical perspective Berry investigated its expected value and variance in \cite{Be:02}, whereas the first mathematically rigorous derivation of the variance was given in \cite{Wi:10}.

In \cite{NPR:19} and \cite{PV:20}, the authors proved central limit theorems, as $k\sto\infty$, for the nodal length of planar Laplacian eigenfunctions, i.e. when $\MM=\R^2$ for $\LL_1^{f_k}(\EE_0(f_k,D))$, in a fixed convex body $D\subset \R^2$, using the so-called fourth moment theorem of \cite{PT:05}. More precisely, showing that the random functional $\LL_1^{f_k}(\EE_0(f_k,D))$ is dominated by the fourth chaotic projection of its Wiener chaos expansion.

At the same time, in \cite{MRW:19}, the authors proved a central limit theorem, as $k\sto\infty$, for the nodal length of Laplacian eigenfunctions on the two-dimensional sphere, i.e. when $\MM=\S^2$ for  $\LL_1^{f_k}(\EE_0(f_k,\S^2))$, using a different idea: 
instead of studying the asymptotic behavior of the entire dominant fourth chaotic component, which is given by a sum of six terms involving the eigenfunctions and their gradients, they proved its asymptotic full correlation with a functional that only depends on the eigenfunction $f_k$ and not on its gradient components. Such functional is the so called \emph{(centered) sample trispectrum} which is defined as the integral of $H_4(f_k)$, where $H_4$ is the fourth Hermite polynomial. This means that the authors of \cite{MRW:19} were able to obtain a much simpler expression for the leading term, making the derivation of a quantitative central limit theorem much more immediate.

Hence a natural question arises: as $k\sto\infty$, is that possible to obtain an asymptotic neater expression also on the plane, that is for $\LL_1^{f_k}(\EE_0(f_k,D))$ when $\MM=\R^2$ and $D\subset\R^2$? Is the fourth chaotic component of the nodal length of planar random wave asymptotically fully correlated with a term that does not depend on the gradient? Here, we will positively answer to these questions, showing that the computations are actually very similar to the ones of \cite{MRW:19}. Indeed, the aim of this short note is not only answering these questions but also highlighting some open ones, that are probably more challenging to address.

In fact, it is important to point out that the asymptotic full correlation of $\LL_1^{f_k}(\EE_0(f_k,\S^2))$ with the (centered) sample trispectrum led to the fact that $\LL_1^{f_k}(\EE_0(f_k,\S^2))$ is also asymptotically fully correlated with the total number of critical points. Indeed, in the paper \cite{CM:19} it is shown that the asymptotic behavior of the total number of critical points is dominated by exactly the same component as the one that dominates in the nodal length, that is the (centered) sample trispectrum. As a consequence, it would be interesting to discover if similar results can be proved in the planar case; heuristics clearly suggest that higher number of critical points would presumably correspond to a higher number of nodal components.

For a threshold parameter $u\ne 0$, asymptotic full correlation of Lipschitz-Killing curvatures and critical values among themselves and with a functional of just the eigenfunction $f_k$ was proved in the works \cite{MW:11, Ro:15b, CM:18, CM:20}, some years before considering the degenerate (and hence more challenging) case $u=0$. Such functional is the so-called \emph{(centered) sample power spectrum}, which is defined as the integral of $H_2(f_k)$, where $H_2$ is the second Hermite polynomial. Moreover, in \cite{MR:19}, it was proved that the correlation between $\LL_1^{f_k}(\EE_0(f_k,\S^2))$ and $\LL_1^{f_k}(\EE_u(f_k,\S^2))$ at any level $u\ne0$ is asymptotically zero, while the partial correlation after controlling for the random norm $\norm{f_k}_{L^2(\S^2)}$ is asymptotically one. In general, it would be interesting to study whether these results hold in the planar case.

\subsection{Main Results}

Let us now get into the notation of \cite{NPR:19,PV:20}. From now on $\MM=\R^2$ and we let
$\Delta$ be the Laplace operator on $\R^{2}$. For $E> 0$, we define 
\begin{equation}\label{rwbk}
B_E(x)=\int_{\mathbb{S}^1}e^{i2\pi^2E\l\theta,x\r}Z(d\theta) \,, \quad x\in \R^2\,,
\end{equation}
where $Z$ is an appropriate Hermitian Gaussian measure on $\mathbb{S}^1$; then $B_E: \R^2 \to \R$
is a Gaussian random field on $\R^2$ such that $\E \[B_E(x)\]=0$ and 
$$
\E\left[B_E(x) B_E(y)\right] = J_0(\sqrt{2\pi^2E}\| x- y\|),\qquad x,y\in \R^2
$$
where $J_0$ denotes the zero-order Bessel function of the first kind (see \cite{Kr:14})
\begin{equation}\label{seriesJ0}
J_0(t) = \sum_{m=0}^{\infty} \frac{(-1)^m}{(m!)^2}\left (\frac{t}{2}\right)^{2m},\qquad t\in \R\,.
\end{equation} 
Moreover, $B_E$ almost surely solves the Helmholtz equation
\begin{equation}\label{Helmholtz}
\Delta B_E+\lambda_E B_E=0\,, \qquad \lambda_E:=2\pi^2E\,,
\end{equation}
so that $B_E=f_k$, in the notation of the previous section.

In this paper, we focus on the nodal length of the random fields
$\{B_E(\cdot)\} ,$ i.e. the boundary length of the excursion set at the level $u=0$ inside a fixed convex body $D\subset \R^2$:
\begin{equation}\label{NL}
\L_{E}:=\LL_1^{B_E}(\EE_0(B_E,D))=\operatorname{length}\{ B_E^{-1}(0)\cap D\}.
\end{equation}%

It is a straightforward application of the Gaussian Kinematic Formula, see \cite{AT:09}, showing that the expectation of the nodal length $\L_E$ satisfies the following relation
\begin{equation}\label{meanlength}
\E[ \L_E] = \area(D)\,\frac{\pi}{\sqrt{2}}\sqrt{E}\,, \quad \forall E>0\,,
\end{equation}
whereas it is more challenging to prove that the variance verifies the asymptotic relation (see \cite{Be:02,Wi:10,NPR:19})
\begin{equation}\label{e:varlength}
\Var(\L_E) =\frac{ \area(D)}{512\pi}\,\log E+o(\log E)\,, \quad E\to\infty\,. 
\end{equation}
Moreover, in \cite{NPR:19}, a CLT is established and, as $E\sto \infty$, it is shown that
\begin{equation*}
\widetilde \L_E :=\frac{\L_E - \E[\L_E]}{\sqrt{\Var(\L_E)}}\,, \mathop{\longrightarrow}^{d} N\,,
\end{equation*}
where $N\sim \mathcal{N}(0,1)$ is a standard Gaussian random variable and $\displaystyle\mathop{\longrightarrow}^{d}$ denotes convergence in distribution. 

In this short note, we will prove the
asymptotic equivalence (in the $L^{2}(\Omega )$ sense) of the nodal length $\L_{E}$
and the (centered) sample trispectrum of $\{B_E\}$, i.e.
\begin{equation}
h_{E ;4}:=\int_{D}H_{4}(B_E(x))dx \,, \label{HE}
\end{equation}
where $H_{4}$ is the fourth-order Hermite polynomial -- we recall that $H_{4}(u)=u^{4}-6u^{2}+3$. 
Now, let us define the following properly rescaled random variables
\begin{equation}
\mathscr{M}_{E}:=-\frac{\sqrt{2\pi^2E}}{96}\,\int_{D}H_{4}(B_E(x))dx=-\frac{\sqrt{2\pi^2E}}{96}\,h_{E ;4} \, .   \label{ME} 
\end{equation}
From \cite[Lemma 8.4]{NPR:19} we know that, as $E \rightarrow \infty $, we have
\begin{equation}\label{eq:trispec var log}
\Var\{ \M_E\} =\frac{\area D}{512\pi}\log E +O(1).
\end{equation}
Looking at \eqref{e:varlength} and \eqref{eq:trispec var log}, it is immediate to note that the variance of $\M_E$ is asymptotically equivalent to the variance of $\L_{E}$, i.e.
\begin{equation*}
\frac{\Var\{ \L_{E}\} }{\Var\{ \M_{E
}\} }=1+o(1)\text{ , as }E \rightarrow \infty\,.
\end{equation*}%
Now, set
\begin{equation}
\widetilde{\M}_{E }:=\frac{\M_E}{\sqrt{\Var(\M_E)}}\,; \label{Mellehat}
\end{equation}
the main result of this note is the following theorem, which is the planar counterpart of \cite[Theorem 1.2]{MRW:19}.

\begin{thm}
\label{maintheorem} As $E \rightarrow \infty ,$ we have that%
\begin{equation}\label{v1}
\mathbb{E}\left[ \{ \widetilde{\L}_{E }-\widetilde{\M}_{E }\} ^{2}\right] =o\left(1\right),
\end{equation}%
and in particular%
\begin{equation*}
\widetilde{\L}_{E }=\widetilde{\M}_{E }+o_{\P}\left(1\right).
\end{equation*}
\end{thm}

The previous result states that the normalized nodal length \eqref{NL} and (centered) sample trispectrum \eqref{ME} are asymptotically equivalent in $L^{2}(\Omega )$, as $E\sto \infty$, and hence in probability and in law.

Now we briefly recall the definition of Wasserstein distance between two random variables $X$
and $Y$, that is the following (see e.g. \cite[Appendix C]{NP:12} for more details)
\begin{equation*}
d_{W}(X,Y):=\sup_{h:\left\Vert h\right\Vert _{Lip}\leq 1}\left\vert
\E[h(X)]-\E[h(Y)]\right\vert\,.
\end{equation*}
It is a well-known fact that convergence in $L^2(\Omega)$ implies convergence in Wasserstein distance, and that
both imply convergence in law.

Let $N\sim \mathcal{N}(0,1)$ be a standard Gaussian random variable; a straightforward consequence of Theorem \ref{maintheorem} is the following quantitative central limit theorem, in Wasserstein distance (see also \cite[Corollary 1.3]{MRW:19}).

\begin{cor}
\label{CLT} As $E \rightarrow \infty ,$ we have that
\begin{equation*}
d_{W}(\widetilde{\L}_{E },N)=o\left(1\right).
\end{equation*}
\end{cor}

The reduction principle in Theorem \ref{maintheorem} also allows to establish Moderate Deviation estimates for the nodal length of planar random waves, see \cite[Remark 1.9]{MRT:20}. The proof of the following result, which is a refinement of the Central Limit Theorem for the sample trispectrum, is analogous to the proof of Lemma 3.1 in \cite{MRT:20} and hence omitted. 
\begin{cor}
Let $\lbrace a_E, E>0 \rbrace$ be any sequence of positive numbers such that, as $E\sto \infty$,
\begin{equation}\label{aE}
a_E \to \infty,\qquad a_E/ (\log E)^{1/14} \to 0.
\end{equation}
Then the sequence of random variables $\lbrace \widetilde{\mathscr M}_E / a_E, E>0\rbrace$ satisfies a Moderate Deviation principle with speed $a_E^2$ and Gaussian rate function $\mathcal I(x):=x^2/2$, $x\in \mathbb R$, i.e., for every Borelian set $B\subset \mathbb R$ it holds that 
\begin{flalign*}
-\inf_{x\in \mathring B} \mathcal I(x) &\le \liminf_{E\sto \infty} \frac{1}{a_E^2} \log \mathbb P \left (\widetilde{\mathscr M}_E / a_E\in B \right )\\
&\le \limsup_{E\sto \infty} \frac{1}{a_E^2} \log \mathbb P \left (\widetilde{\mathscr M}_E / a_E\in B \right ) \le -\inf_{x\in \bar B} \mathcal I(x),
\end{flalign*} 
where $\mathring B$ (resp. $\bar B$) denotes the interior (resp. the closure) of $B$. 
\end{cor} 
As for the proof of Theorem 1.7 in \cite{MRT:20}, the two sequences of random variables $\lbrace \widetilde{\mathscr M}_E / a_E, E>0\rbrace$ and $\lbrace \widetilde{\mathscr L}_E / a_E, E>0\rbrace$
being exponentially equivalent \cite[Definition 4.2.10]{DZ:98} as soon as $a_E$ goes to infinity sufficiently slowly (according to both \eqref{v1} and \eqref{aE}),  Moderate Deviation estimates can be deduced for $\{\widetilde{\mathscr L}_E / a_E, E>0\}$ with speed $a_{E}^2$ and Gaussian rate function $\mathcal I$.


\section{\label{ProofTheorems} Proofs}
The proofs of the two results Theorem \ref{maintheorem} and Corollary \ref{CLT} are very similar to the ones in \cite{MRW:19} and are strongly based on various results already proved in \cite{NPR:19}.
\vspace{-.2cm}
\paragraph{\label{L2Discussion}The Wiener Chaos Decomposition of the Nodal Length}
In \cite{NPR:19}, the chaotic expansion of the nodal length in a fixed convex body $D\subset \R^2$ is established:
\begin{equation}\label{expL}
\begin{split}
\L_E = \sum_{q=0}^{\infty} \L_E[2q] 
=&\sqrt{2\pi^2E} \sum_{q=0}^{\infty} \sum_{u=0}^{q} \sum_{m=0}^{u} \beta_{2q - 2u} \alpha_{2m, 2u - 2m}\times \cr
 &\times \int_{D} H_{2q-2u}(B_E(x)) H_{2m}(\widetilde \partial_1 B_E(x)) H_{2u-2m}(\widetilde \partial_2 B_E(x))\,dx,
\end{split}
\end{equation}
where the series converges in $L^2(\Omega)$ and $\lbrace \beta_{2n}\rbrace_{n\ge 0}$ is defined in equation (3.50) of \cite{NPR:19}, while $\lbrace \alpha_{2n,2m}\rbrace_{n,m\ge 0}$ is the sequence of chaotic coeffients of the Euclidean norm in $\R^2$ $\| \cdot \|$ appearing in \cite[Lemma 3.5]{MPRW:16}.
Once the chaotic expansions were established, the authors of \cite{NPR:19} proved that, as $E\sto \infty$, 
\begin{equation*}
\frac{\L_E - \E[\L_E]}{\sqrt{\Var(\L_E)}} = \frac{\L_E[4]}{\sqrt{\Var(\L_E[4])}} + o_{\P}(1)\,, 
\end{equation*}
noting that $\lim_{E\sto\infty}\Var \L_E/ \Var \L_E[4]=1$. In particular, the fourth chaotic component of $\L_E$ is given by
\begin{equation}\label{L_4}
\L_E[4](D)=\frac{\sqrt{2\pi^2\,E}}{128}\{8\,a_{1,E}-a_{2,E}-a_{3,E}-2\,a_{4,E}-8\,a_{5,E}-8\,a_{6,E}\}\,,
\end{equation}
where 
\begin{equation}\label{a's}
\begin{split}
a_{1,E}&:=\int_{D} H_4(B_E(x))dx\,,\quad a_{2,E}:=\int_{D}  H_4(\widetilde{\partial}_1 B_E(x))dx\,,\quad
a_{3,E}:=\int_{D} H_4(\widetilde{\partial}_2 B_E(x))dx\,,\\
a_{4,E} &:= \int_{D}  H_2(\widetilde{\partial}_1 B_E(x))
H_2(\widetilde{\partial}_2 B_E(x))dx\,,\\
a_{5,E} &:= \int_{D}  H_2(B_E(x))
  H_2(\widetilde{\partial}_1 B_E(x))dx\,,\quad
a_{6,E} := \int_{D}  H_2(B_E(x))H_2(\widetilde{\partial}_2 B_E(x))dx\,.
\end{split}
\end{equation}
As proved in \cite[Proposition 6.1]{NPR:19}, its variance satisfies 
\begin{equation}\label{var4L}
\begin{split}
{\Var}(\L_E[4]) &= \frac{\pi^2E}{8192}\,{\Var}\left (8a_{1,E}-a_{2,E}-a_{3,E}-2a_{4,E}-8a_{5,E}-8a_{6,E}
\right )\\
&= \frac{{\area}(D)}{512\pi}\,\log E+O(1)\, \quad \text{ as } \,\,E\to \infty\,. 
\end{split}
\end{equation}



\paragraph{Proof of Theorem \protect\ref{maintheorem}}

To establish Theorem \ref{maintheorem}, it is sufficient to show that,
as $E \rightarrow \infty$,
\begin{equation*}
\Corr\( \L_{E},\M_E\) \to1 \,.
\end{equation*}
We have 
\begin{flalign*}
&\Corr\( \L_{E},\M_E\)=\frac{\Cov\( \L_{E},\M_E\)}{\sqrt{\Var (\L_{E})\Var(\M_E)}}\\
&=\frac{\frac{\log E}{512\pi}+O(1)}{\sqrt{\(\frac{\log E}{512\pi}+o(\log E)\)\(\frac{\log E}{512\pi}+g_E\)}} \qquad \text{ with } \,\, \abs{g_E}\le c\,, \text{ a constant independent of } E \\
&=1+o(1)\,.
\end{flalign*}
Indeed, since $\M_E$ is an element of the fourth Wiener chaos, by orthogonality we have that
\begin{flalign*}
&\Cov\( \L_{E},\M_E\)=\Cov\( \sum_{q\ge0}\L_{E}[2q],\M_E\)=\Cov\(\L_{E}[4],\M_E\)\\
&=\Cov\(\frac{\sqrt{2\pi^2\,E}}{128}\{8\,a_{1,E}-a_{2,E}-a_{3,E}-2\,a_{4,E}-8\,a_{5,E}-8\,a_{6,E}\},- \frac{\sqrt{2\pi^2E}}{96}\,a_{1,E}\)\\
&=\frac{2\pi^2\,E}{(128)(96)}\Cov\(\{8\,a_{1,E}-a_{2,E}-a_{3,E}-2\,a_{4,E}-8\,a_{5,E}-8\,a_{6,E}\},-a_{1,E}\)\\
&=\frac{2\pi^2\,E}{(128)(96)}\Cov\(\{-8\,a_{1,E}+a_{2,E}+a_{3,E}+2\,a_{4,E}+8\,a_{5,E}+8\,a_{6,E}\},a_{1,E}\)\\
&=\frac{2\pi^2\,E}{(128)(96)}\,\bigg[-8\,\Var(a_{1,E})+\Cov\(a_{1,E},a_{2,E}\)+\Cov\(a_{1,E},a_{3,E}\)\\
&\qquad\qquad\qquad +2\Cov\(a_{1,E},a_{4,E}\)+8\,\Cov\(a_{1,E},a_{5,E}\)+8\,\Cov\(a_{1,E},a_{6,E}\)\bigg]\,.
\end{flalign*}
After these simple steps, the fact that $\Cov\( \L_{E},\M_E\)=\frac{\log E}{512\pi}+O(1)$ follows straightforwardly using \cite[Lemma 8.4]{NPR:19}.
\endmproof

%
%

\paragraph{Proof of the Central Limit Theorem (Corollary \protect\ref{CLT})}

The sequence $\{ \widetilde{\M}_{E }\} $ is standardized and belongs to the $4$th Wiener chaos, so that we can apply the fourth moment theorem by Nourdin and Peccati, see \cite[Theorem 5.2.6]{NP:12}, to have
\begin{equation*}
d_{W}(\widetilde{\M}_{E },\mathcal{N}(0,1))\leq \sqrt{\frac{1}{%
2\pi }\,\cum_4\(\widetilde{\M}_{E }\) }\,.
\end{equation*}%
Now, using the diagram formula \cite[Proposition 4.15]{MP:11}, we get\footnote{note that here the computations follow the same ideas of \cite[Lemma 8.1]{NPR:19}, however, we are making the details more evident.} that
\begin{flalign*}
&\cum_4\(\widetilde{\M}_{E }\)=\frac{\pi^4E^2}{81(\area D)^2\log^2 E}\times\\
&\times\,\int_{D^4}\cum_4\(H_4\(B_E(x)\),H_4\(B_E(y)\),H_4\(B_E(z)\),H_4\(B_E(w)\)\)\, dx\,dy\,dz\,dw\\
&=\frac{\pi^4E^2}{81(\area D)^2\log^2 E}\,\int_{D^4} \bigg[6^5\,r^E(x-y)^2\,r^E(y-z)^2\,r^E(z-w)^2\,r^E(w-x)^2+\\
&\qquad\qquad+\(6 \cdot 4^4\)\,r^E(x-y)^3\,r^E(y-z)\,r^E(z-w)^3\,r^E(w-x)\bigg]\, dx\,dy\,dz\,dw\\
&=\frac{96\pi^4E^2}{(\area D)^2\log^2 E}\,\int_{D^4} r^E(x-y)^2\,r^E(y-z)^2\,r^E(z-w)^2\,r^E(w-x)^2\, dx\,dy\,dz\,dw\\
&\quad+\frac{512\pi^4E^2}{27(\area D)^2\log^2 E}\,\int_{D^4}\,r^E(x-y)^3\,r^E(y-z)\,r^E(z-w)^3\,r^E(w-x)\, dx\,dy\,dz\,dw\\
&=\frac{96\pi^4}{(\area D)^2E^2\log^2 E}\,\int_{ (\sqrt E D)^4} r^1(x-y)^2\,r^1(y-z)^2\,r^1(z-w)^2\,r^1(w-x)^2\, dx\,dy\,dz\,dw\\
&\quad+\frac{512\pi^4}{27(\area D)^2E^2\log^2 E}\,\int_{ (\sqrt E D)^4}\,r^1(x-y)^3\,r^1(y-z)\,r^1(z-w)^3\,r^1(w-x)\, dx\,dy\,dz\,dw\\
&:=\mathcal I_1^E+ \mathcal I_2^E\,.
\end{flalign*}
For the first summand, we have 
\begin{flalign*}
\mathcal I_1^E&=\frac{96\pi^4}{(\area D)^2E^2\log^2 E}\,\int_{ (\sqrt E D)^4} J_0\(2\pi\norm{x_1-x_2}\)^2\,J_0\(2\pi\norm{x_2-x_3}\)^2\\
&\qquad \qquad\qquad \qquad\,J_0\(2\pi\norm{x_3-x_4}\)^2\,J_0\(2\pi\norm{x_4-x_1}\)^2\, dx_1\,dx_2\,dx_3\,dx_4\\
&\le \frac{2\cdot96\,\pi^4}{(\area D)^2E^2\log^2 E}\,\int_{ (\sqrt E D)^4} J_0\(2\pi\norm{x_1-x_2}\)^4\\
&\qquad \qquad\qquad \qquad\,\,J_0\(2\pi\norm{x_3-x_4}\)^2\,J_0\(2\pi\norm{x_4-x_1}\)^2\, dx_1\,dx_2\,dx_3\,dx_4\,.
\end{flalign*}
Now, we make the following change of variable
\begin{equation}
\begin{split}
x_{2,1}=x_{1,1}+\phi \cos\theta \quad x_{2,2}=x_{1,2}+\phi \sin\theta\,, \quad \phi \in (0,\sqrt{E}\diam D]\,, \, \theta\in[0,2\pi] \\
x_{3,1}=x_{4,1}+\rho \cos\alpha \quad x_{3,2}=x_{4,2}+\rho \sin\alpha\,, \quad \rho \in (0,\sqrt{E}\diam D]\,, \, \alpha\in[0,2\pi]\\
x_{1,1}=x_{4,1}+\psi \cos\beta \quad x_{1,2}=x_{4,2}+\psi \sin\beta\,, \quad \psi \in (0,\sqrt{E}\diam D]\,, \, \beta\in[0,2\pi]
\end{split}
\end{equation}
to have
\begin{flalign*}
\mathcal I_1^E&\le \frac{2\cdot96\,\pi^4}{(\area D)^2E^2\log^2 E}\,\int_{\sqrt E D}\,dx_4 \, \int_{[0,2\pi]^3} d\theta\,d\alpha\,d\beta \int_{(0,\sqrt{E}\diam D]^3}\,\\
&\qquad \qquad J_0\(2\pi\phi\)^4\,J_0\(2\pi\rho\)^2\,J_0\(2\pi\psi\)^2\, \phi \,d\phi \,\psi \,d\psi\,\rho\, d\rho\\
&= \frac{2\cdot96\,\pi^4}{(\area D)^2E^2\log^2 E}\,\int_{\sqrt E D}\,dx_4 \, \int_{[0,2\pi]^3} d\theta\,d\alpha\,d\beta \int_{(1,\sqrt{E}\diam D]^3}\,\\
&\qquad \qquad J_0\(2\pi\phi\)^4\,J_0\(2\pi\rho\)^2\,J_0\(2\pi\psi\)^2\, \phi \,d\phi \,\psi \,d\psi\,\rho\, d\rho+O\(\frac1{\log^2 E}\)\\
&= \frac{2^4\cdot96\,\pi^7}{\area D\,E\,\log^2 E} \int_{(1,\sqrt{E}\diam D]^3}\, \frac1\phi \,d\phi  \,d\psi\, d\rho+O\(\frac1{\log^2 E}\)=O\(\frac1{\log E}\)\,,
\end{flalign*}
where we used the fact that, as $r\sto\infty$, 
$$
J_0(2\pi r)=\frac{1}{\pi\sqrt{r}}\cos\(2\pi r-\frac \pi 4\)+O\(\frac{1}{r^{3/2}}\)\,.
$$
Analogously, for the second summand, we have that
\begin{flalign*}
\mathcal I_2^E&=\frac{512\pi^4}{27(\area D)^2E^2\log^2 E}\,\int_{ (\sqrt E D)^4} J_0\(2\pi\norm{x_1-x_2}\)^3\,J_0\(2\pi\norm{x_2-x_3}\)\notag\\
&\qquad \qquad\qquad \qquad\,J_0\(2\pi\norm{x_3-x_4}\)^3\,J_0\(2\pi\norm{x_4-x_1}\)\, dx_1\,dx_2\,dx_3\,dx_4\notag\\
&\le \frac{2\cdot512\,\pi^4}{27(\area D)^2E^2\log^2 E}\,\int_{ (\sqrt E D)^4} J_0\(2\pi\norm{x_1-x_2}\)^4\notag\\
&\qquad \qquad\qquad \qquad\,\abs{J_0\(2\pi\norm{x_3-x_4}\)}^3\,\abs{J_0\(2\pi\norm{x_4-x_1}\)}\, dx_1\,dx_2\,dx_3\,dx_4\notag\\
&=\frac{2\cdot512\,\pi^4}{27(\area D)^2E^2\log^2 E}\,\int_{\sqrt E D}\,dx_4 \, \int_{[0,2\pi]^3} d\theta\,d\alpha\,d\beta \int_{(0,\sqrt{E}\diam D]^3}\,\notag\\
&\qquad \qquad J_0\(2\pi\phi\)^4\,\abs{J_0\(2\pi\rho\)}^3\,\abs{J_0\(2\pi\psi\)}\, \phi \,d\phi \,\psi \,d\psi\,\rho\, d\rho\notag\\
&=\frac{2\cdot512\,\pi^4}{27(\area D)^2E^2\log^2 E}\,\int_{\sqrt E D}\,dx_4 \, \int_{[0,2\pi]^3} d\theta\,d\alpha\,d\beta \int_{(1,\sqrt{E}\diam D]^3}\,\notag\\
&\qquad \qquad J_0\(2\pi\phi\)^4\,\abs{J_0\(2\pi\rho\)}^3\,\abs{J_0\(2\pi\psi\)}\, \phi \,d\phi \,\psi \,d\psi\,\rho\, d\rho+O\(\frac1{\log^2 E}\)\notag\\
&=\frac{2^4\cdot512\,\pi^7}{27\area D\,E\log^2 E}\, \int_{(1,\sqrt{E}\diam D]^3}\, \frac{1}{\phi}\, \frac{1}{\rho^{1/2}}\,\psi^{1/2} \,d\phi  \,d\psi\, d\rho+O\(\frac1{\log^2 E}\)=O\(\frac1{\log E}\).
\end{flalign*}
Consequently, we just proved that
\begin{equation}\label{Ogrande}
\cum_4\(\widetilde{\M}_{E }\)= \mathcal I_1^E+\mathcal I_2^E=O\(\frac1{\log E}\)
\end{equation}
Thanks to the triangle inequality valid for the Wasserstein distance $d_{W}$ (see \cite[Appendix C]{NP:12}), we have that
\begin{flalign*}
d_{W}(\widetilde{\L}_{E },\mathcal{N}(0,1))&\leq d_{W}(\widetilde{\M}_{E },\mathcal{N}(0,1))+\sqrt{\mathbb{E}\left[ \(\widetilde{%
\mathcal{L}}_{E }-\widetilde{\M}_{E }\)^{2}\right]}\\
&=O\left(\frac{1}{\sqrt{\log E }}\right)+o(1)=o(1)\,, 
\end{flalign*}
as $E\sto\infty$, where the first equality follows from \eqref{Ogrande} and Theorem \ref{maintheorem}. The proof is hence concluded.
\endmproof

\paragraph*{Acknowledgements} 
The author would like to thank Domenico Marinucci and Maurizia Rossi for useful discussions. This note has been written when the author was PostDoc at the Department of Mathematics, University of Rome Tor Vergata, and she acknowledges also the MIUR Excellence Department Project, CUP E83C18000100006.

\end{document}